% Submitted to the JSL on July 9, 1996
%
% Corrections, suggested by the referee, made on December 23, 1996, 
% and resubmitted as the final version on December 23, 1996.
% 
%\input fonts
\font\fifteenrm=cmr10 scaled\magstep2 % this is all really just 14.4pt
\font\fifteeni=cmmi10 scaled\magstep2
\font\fifteensy=cmsy10 scaled\magstep2
\font\fifteenbf=cmbx10 scaled\magstep2
\font\fifteentt=cmtt10 scaled\magstep2
\font\fifteenit=cmti10 scaled\magstep2
\font\fifteensl=cmsl10 scaled\magstep2
\font\fifteenam=msam10 scaled\magstep2
\font\fifteenbm=msbm10 scaled\magstep2
\font\fifteenex=cmex10 scaled\magstep2
\font\fifteensc=cmcsc10 scaled\magstep2 
\font\twelverm=cmr10 at 12pt
\font\twelvei=cmmi10 at 12pt
\font\twelvesy=cmsy10 at 12pt
\font\twelvebf=cmbx10 at 12pt
\font\twelvett=cmtt10 at 12pt
\font\twelveit=cmti10 at 12pt
\font\twelvesl=cmsl10 at 12pt
\font\twelveam=msam10 at 12pt
\font\twelvebm=msbm10 at 12pt
\font\twelveex=cmex10 at 12pt
\font\twelvesc=cmcsc10 at 12pt
\font\elevenrm=cmr10 scaled\magstephalf % this is really 10.95 pt
\font\eleveni=cmmi10 scaled\magstephalf
\font\elevensy=cmsy10 scaled\magstephalf
\font\elevenbf=cmbx10 scaled\magstephalf
\font\eleventt=cmtt10 scaled\magstephalf
\font\elevenit=cmti10 scaled\magstephalf
\font\elevensl=cmsl10 scaled\magstephalf
\font\elevenam=msam10 scaled\magstephalf
\font\elevenbm=msbm10 scaled\magstephalf
\font\elevenex=cmex10 scaled\magstephalf
\font\elevensc=cmcsc10 scaled\magstephalf
\font\tenrm=cmr10
\font\teni=cmmi10
\font\tensy=cmsy10
\font\tenbf=cmbx10
\font\tentt=cmtt10
\font\tenit=cmti10
\font\tensl=cmsl10
\font\tenam=msam10
\font\tenbm=msbm10
\font\tenex=cmex10
\font\tensc=cmcsc10
\font\ninerm=cmr9
\font\ninei=cmmi9
\font\ninesy=cmsy9
\font\ninebf=cmbx9
\font\ninett=cmtt9
\font\nineit=cmti9
\font\ninesl=cmsl9
\font\nineam=msam9
\font\ninebm=msbm9
\font\nineex=cmex9
\font\ninesc=cmcsc9
\font\eightrm=cmr8
\font\eighti=cmmi8
\font\eightsy=cmsy8
\font\eightbf=cmbx8
\font\eighttt=cmtt8
\font\eightit=cmti8
\font\eightsl=cmsl8
\font\eightam=msam8
\font\eightbm=msbm8
\font\eightex=cmex8
\font\eightsc=cmcsc8
\font\sevenrm=cmr7
\font\seveni=cmmi7
\font\sevensy=cmsy7
\font\sevenbf=cmbx7

\font\sevenam=msam7
\font\sevenbm=msbm7

\font\sixrm=cmr6
\font\sixi=cmmi6
\font\sixsy=cmsy6

\font\sixam=msam6
\font\sixbm=msbm6

\font\fiverm=cmr5
\font\fivei=cmmi5
\font\fivesy=cmsy5

\font\fiveam=msam5
\font\fivebm=msbm5

\font\fourrm=cmr5 at 4pt
\font\fouri=cmmi5 at 4pt
\font\foursy=cmsy5 at 4pt

\font\fouram=msam5 at 4pt
\font\fourbm=msbm5 at 4pt

\skewchar\twelvei='177 \skewchar\eleveni='177\skewchar\teni='177
\skewchar\ninei='177 \skewchar\eighti='177\skewchar\seveni='177 
\skewchar\sixi='177 \skewchar\fivei='177 \skewchar\fouri='177
\skewchar\twelvesy='60 \skewchar\elevensy='60 \skewchar\tensy='60
\skewchar\ninesy='60 \skewchar\eightsy='60 \skewchar\sevensy='60 
\skewchar\sixsy='60 \skewchar\fivesy='60 \skewchar\foursy='60
\newfam\itfam
\newfam\slfam
\newfam\bffam
\newfam\ttfam
\newfam\scfam
\newfam\amfam
\newfam\bmfam
\def\eightbig#1{{\hbox{$\left#1\vbox to 6.5pt{}\voidright $}}}
\def\eightBig#1{{\hbox{$\left#1\vbox to 7.5pt{}\voidright $}}}
\def\eightbigg#1{{\hbox{$\left#1\vbox to 10pt{}\voidright $}}}
\def\eightBigg#1{{\hbox{$\left#1\vbox to 13pt{}\voidright $}}}
\def\ninebig#1{{\hbox{$\left#1\vbox to 7.5pt{}\voidright $}}}
\def\nineBig#1{{\hbox{$\left#1\vbox to 8.5pt{}\voidright $}}}
\def\ninebigg#1{{\hbox{$\left#1\vbox to 11.5pt{}\voidright $}}}
\def\nineBigg#1{{\hbox{$\left#1\vbox to 14.5pt{}\voidright $}}}
\def\tenbig#1{{\hbox{$\left#1\vbox to 8.5pt{}\voidright $}}}
\def\tenBig#1{{\hbox{$\left#1\vbox to 9.5pt{}\voidright $}}}
\def\tenbigg#1{{\hbox{$\left#1\vbox to 12.5pt{}\voidright $}}}
\def\tenBigg#1{{\hbox{$\left#1\vbox to 16pt{}\voidright $}}}
\def\elevenbig#1{{\hbox{$\left#1\vbox to 9pt{}\voidright $}}}
\def\elevenBig#1{{\hbox{$\left#1\vbox to 10.5pt{}\voidright $}}}
\def\elevenbigg#1{{\hbox{$\left#1\vbox to 14pt{}\voidright $}}}
\def\elevenBigg#1{{\hbox{$\left#1\vbox to 17.5pt{}\voidright $}}}
\def\twelvebig#1{{\hbox{$\left#1\vbox to 10pt{}\voidright $}}}
\def\twelveBig#1{{\hbox{$\left#1\vbox to 11pt{}\voidright $}}}
\def\twelvebigg#1{{\hbox{$\left#1\vbox to 15pt{}\voidright $}}}
\def\twelveBigg#1{{\hbox{$\left#1\vbox to 19pt{}\voidright $}}}
\def\fifteenbig#1{{\hbox{$\left#1\vbox to 12pt{}\voidright $}}}
\def\fifteenBig#1{{\hbox{$\left#1\vbox to 13.5pt{}\voidright $}}}
\def\fifteenbigg#1{{\hbox{$\left#1\vbox to 18pt{}\voidright $}}}
\def\fifteenBigg#1{{\hbox{$\left#1\vbox to 23pt{}\voidright $}}}
\def\voidright{\right.\nulldelimiterspace=0pt \mathsurround=0pt }
\def\fifteenpoint{
  \textfont0=\fifteenrm \scriptfont0=\twelverm \scriptscriptfont0=\tenrm
  \def\rm{\fam0 \fifteenrm}%
  \textfont1=\fifteeni \scriptfont1=\twelvei \scriptscriptfont1=\teni
  \textfont2=\fifteensy \scriptfont2=\twelvesy \scriptscriptfont2=\tensy
  \textfont3=\fifteenex \scriptfont3=\fifteenex \scriptscriptfont3=\fifteenex
  \def\it{\fam\itfam\fifteenit}\textfont\itfam=\fifteenit
  \def\sl{\fam\slfam\fifteensl}\textfont\slfam=\fifteensl
  \def\bf{\fam\bffam\fifteenbf}\textfont\bffam=\fifteenbf 
    \scriptfont\bffam=\twelvebf\scriptscriptfont\bffam=\tenbf
  \def\tt{\fam\ttfam\fifteentt}\textfont\ttfam=\fifteentt
  \def\sc{\fam\scfam\fifteensc}\textfont\scfam=\fifteensc
  \def\am{\fam\amfam\fifteenam}\textfont\amfam=\fifteenam
    \scriptfont\amfam=\twelveam\scriptscriptfont\amfam=\tenam
  \def\bm{\fam\bmfam\fifteenbm}\textfont\bmfam=\fifteenbm
    \scriptfont\bmfam=\twelvebm\scriptscriptfont\bmfam=\tenbm
  \baselineskip=21pt \rm
  \let\big=\fifteenbig\let\Big=\fifteenBig\let\bigg=\fifteenbigg
  \let\Bigg=\fifteenBigg}
\def\twelvepoint{
  \textfont0=\twelverm \scriptfont0=\ninerm \scriptscriptfont0=\sevenrm
  \def\rm{\fam0 \twelverm}%
  \textfont1=\twelvei \scriptfont1=\ninei \scriptscriptfont1=\seveni
  \textfont2=\twelvesy \scriptfont2=\ninesy \scriptscriptfont2=\sevensy
  \textfont3=\twelveex \scriptfont3=\twelveex \scriptscriptfont3=\twelveex
  \def\it{\fam\itfam\twelveit}\textfont\itfam=\twelveit
  \def\sl{\fam\slfam\twelvesl}\textfont\slfam=\twelvesl
  \def\bf{\fam\bffam\twelvebf}\textfont\bffam=\twelvebf 
    \scriptfont\bffam=\ninebf\scriptscriptfont\bffam=\sevenbf
  \def\tt{\fam\ttfam\twelvett}\textfont\ttfam=\twelvett
  \def\sc{\fam\scfam\twelvesc}\textfont\scfam=\twelvesc
  \def\am{\fam\amfam\twelveam}\textfont\amfam=\twelveam
    \scriptfont\amfam=\nineam\scriptscriptfont\amfam=\sevenam
  \def\bm{\fam\bmfam\twelvebm}\textfont\bmfam=\twelvebm
    \scriptfont\bmfam=\ninebm\scriptscriptfont\bmfam=\sevenbm
  \baselineskip=17.8pt \rm 
  \def\looselineskip{\baselineskip=18.5pt plus 1.8pt}%
  \def\tightlineskip{\baselineskip=16.5pt}%
  \def\verytightlineskip{\baselineskip=15pt}%
  \let\big=\twelvebig\let\Big=\twelveBig\let\bigg=\twelvebigg
  \let\Bigg=\twelveBigg  }
\def\elevenpoint{
  \textfont0=\elevenrm \scriptfont0=\ninerm \scriptscriptfont0=\sixrm
  \def\rm{\fam0 \elevenrm}%
  \textfont1=\eleveni \scriptfont1=\ninei \scriptscriptfont1=\sixi
  \textfont2=\elevensy \scriptfont2=\ninesy \scriptfont2=\sixsy 
  \textfont3=\elevenex \scriptfont3=\elevenex \scriptfont3=\elevenex
  \def\it{\fam\itfam\elevenit}\textfont\itfam=\elevenit
  \def\sl{\fam\slfam\elevensl}\textfont\slfam=\elevensl
  \def\bf{\fam\bffam\elevenbf}\textfont\bffam=\elevenbf
  \def\tt{\fam\ttfam\eleventt}\textfont\ttfam=\eleventt
  \def\sc{\fam\scfam\elevensc}\textfont\scfam=\elevensc
  \def\am{\fam\amfam\elevenam}\textfont\amfam=\elevenam
    \scriptfont\amfam=\nineam\scriptscriptfont\amfam=\sixam
  \def\bm{\fam\bmfam\elevenbm}\textfont\bmfam=\elevenbm
    \scriptfont\bmfam=\ninebm\scriptscriptfont\bmfam=\sixbm
  \baselineskip=15.1pt \rm
  \def\looselineskip{\baselineskip=16pt plus 1.5pt}%
  \def\tightlineskip{\baselineskip=14pt}%
  \def\verytightlineskip{\baselineskip=13pt}%
  \let\big=\elevenbig\let\Big=\elevenBig\let\bigg=\elevenbigg
  \let\Bigg=\elevenBigg  }
\def\tenpoint{
  \textfont0=\tenrm \scriptfont0=\eightrm \scriptscriptfont0=\fiverm
  \def\rm{\fam0 \tenrm}%
  \textfont1=\teni \scriptfont1=\eighti \scriptscriptfont1=\fivei
  \textfont2=\tensy \scriptfont2=\eightsy \scriptfont2=\fivesy 
  \textfont3=\tenex \scriptfont3=\tenex \scriptfont3=\tenex
  \def\it{\fam\itfam\tenit}\textfont\itfam=\tenit
  \def\sl{\fam\slfam\tensl}\textfont\slfam=\tensl
  \def\bf{\fam\bffam\tenbf}\textfont\bffam=\tenbf
  \def\tt{\fam\ttfam\tentt}\textfont\ttfam=\tentt
  \def\sc{\fam\scfam\tensc}\textfont\scfam=\tensc
  \def\am{\fam\amfam\tenam}\textfont\amfam=\tenam
    \scriptfont\amfam=\eightam \scriptscriptfont\amfam=\fiveam
  \def\bm{\fam\bmfam\tenbm}\textfont\bmfam=\tenbm
    \scriptfont\bmfam=\eightbm \scriptscriptfont\bmfam=\fivebm
  \baselineskip=14pt\rm
  \def\looselineskip{\baselineskip=14.8pt plus1.5pt}
  \def\tightlineskip{\baselineskip=13.6pt}%
  \def\verytightlineskip{\baselineskip=13pt}%
  \let\big=\tenbig\let\Big=\tenBig\let\bigg=\tenbigg\let\Bigg=\tenBigg  }
\def\ninepoint{
  \textfont0=\ninerm \scriptfont0=\sevenrm \scriptscriptfont0=\fourrm
  \def\rm{\fam0 \ninerm}%
  \textfont1=\ninei \scriptfont1=\seveni \scriptscriptfont1=\fouri
  \textfont2=\ninesy \scriptfont2=\sevensy \scriptfont2=\foursy 
  \textfont3=\nineex \scriptfont3=\nineex \scriptfont3=\nineex
  \def\it{\fam\itfam\nineit}\textfont\itfam=\nineit
  \def\sl{\fam\slfam\ninesl}\textfont\slfam=\ninesl
  \def\bf{\fam\bffam\ninebf}\textfont\bffam=\ninebf
  \def\tt{\fam\ttfam\ninett}\textfont\ttfam=\ninett
  \def\sc{\fam\scfam\ninesc}\textfont\scfam=\ninesc
  \def\am{\fam\amfam\nineam}\textfont\amfam=\nineam
    \scriptfont\amfam=\nineam\scriptscriptfont\amfam=\fouram
  \def\bm{\fam\bmfam\ninebm}\textfont\bmfam=\ninebm
    \scriptfont\bmfam=\ninebm\scriptscriptfont\bmfam=\fourbm
  \baselineskip=12.6pt\rm
  \def\tightlineskip{\baselineskip=11.5pt}
  \let\big=\ninebig\let\Big=\nineBig\let\bigg=\ninebigg
  \let\Bigg=\nineBigg  }
\def\eightpoint{
  \textfont0=\eightrm \scriptfont0=\fiverm \scriptscriptfont0=\fourrm
  \def\rm{\fam0 \eightrm}%
  \textfont1=\eighti \scriptfont1=\fivei \scriptscriptfont1=\fouri
  \textfont2=\eightsy \scriptfont2=\fivesy \scriptfont2=\foursy 
  \textfont3=\eightex \scriptfont3=\eightex \scriptfont3=\eightex
  \def\it{\fam\itfam\eightit}\textfont\itfam=\eightit
  \def\sl{\fam\slfam\eightsl}\textfont\slfam=\eightsl
  \def\bf{\fam\bffam\eightbf}\textfont\bffam=\eightbf
  \def\tt{\fam\ttfam\eighttt}\textfont\ttfam=\eighttt
  \def\sc{\fam\scfam\eightsc}\textfont\scfam=\eightsc
  \def\am{\fam\amfam\eightam}\textfont\amfam=\eightam
    \scriptfont\amfam=\eightam\scriptscriptfont\amfam=\fouram
  \def\bm{\fam\bmfam\eightbm}\textfont\bmfam=\eightbm
    \scriptfont\bmfam=\eightbm\scriptscriptfont\bmfam=\fourbm
  \baselineskip=11.2pt \rm
  \let\big=\eightbig\let\Big=\eightBig\let\bigg=\eightbigg
  \let\Bigg=\eightBigg  }

\twelvepoint
\nopagenumbers
\hsize=6in\vsize=8.8in

\parskip=1pt plus 1pt

\newif\ifSpecialhead\Specialheadfalse
\newbox\specialheadbox

\def\specialhead #1\par{\Specialheadtrue\setbox\specialheadbox=\hbox{#1}}
\headline={{\ifSpecialhead\box\specialheadbox\global\Specialheadfalse\else
     \ifnum\pageno<0{\hfill\quad{\twelvebf\folio}}%
     \else\ifnum\pageno<2\hfill
     \else\hfill\twelvepoint\sc\firstmark\quad{\twelvebf\folio}\fi\fi\fi}}

\def\title#1\par{\bigskip{\def\cr{\par\center}\center\fifteenbf #1\par}\medskip}
\def\subtitle#1\par{\centerline{\fifteenrm #1}\medskip}
\def\author#1\par{\medskip{\def\cr{\par\center\twelvesc}\fifteensc\center#1\par}}
\def\center#1\par{\hfil #1\hfil\par}
\def\abstract.#1\par{\message{Abstract.}%
                    \medskip{\narrower\narrower\tenpoint\tightlineskip
                        \noindent{\bf Abstract.}#1\par}\medskip\noindent}
\def\bigabstract.#1\par{\message{Abstract.}%
                         \medskip{\narrower\narrower\tightlineskip
                         \noindent{\bf Abstract. }#1\par}\medskip\noindent}
\def\acknowledgement#1\par{\footnote{}{#1}}
\def\sectionskip{\Goodbreak\vskip 25pt plus 15pt minus 5pt}
\def\secnumber{\ifquiet
               \else\ifNoSections
                    \else\sectionsymbol\the\secno\quad\fi\fi}
\def\section#1\par{ \NoSectionsfalse\par\sectionskip\proofdepth=0\claimno=0
 \ifquiet\else\advance\secno by1\fi\toks0={#1}
 \immediate\write16{\ifquiet\else Section \the\secno\space\fi
                    \the\toks0}%
 \mark{\secnumber #1}%
 {\fifteenpoint\bf\noindent\secnumber #1}\nobreak\bigskip\quietoff
 \nobreak\noindent}

\def\QUIET{\QUIETtrue\quiettrue}

\def\quietoff{\ifQUIET\else\quietfalse\fi}
\newif\ifquiet
\newif\ifQUIET
\newif\ifNoSections
\newcount\claimtype
\newcount\secno
\newcount\claimno
\newcount\subclaimno
\newcount\subsubclaimno
\newcount\subsubsubclaimno
\newcount\proofdepth
\def\subclaimnumber{\ifquiet\else\ifcase\subclaimno\or A\or B\or C\or D\or E\or
     F\or G\or H\or I\or J\or K\or L\or M\or N\or O\or P\fi\fi}
\def\subsubclaimnumber{\ifquiet\else\ifcase\subsubclaimno\or i\or ii\or iii\or 
   iv\or v\or vi\or vii\or viii\or ix\or x\or xi\or xii\or xiii\or xiv\fi\fi}
\def\subsubsubclaimnumber{\ifquiet\else\ifcase\subsubsubclaimno\or a\or b\or 
   c\or d\or e\or f\or g\or \or h\or i\or j\or k\or l\or m\or n\or o\fi\fi}
\def\claimtag{\ifquiet\else
  \ifNoSections
    \ifcase\proofdepth\the\claimno%
    \or\the\claimno.\subclaimnumber
    \or\the\claimno.\subclaimnumber.\subsubclaimnumber
    \or\the\claimno.\subclaimnumber.\subsubclaimnumber
                                                .\subsubsubclaimnumber\fi
  \else
    \ifcase\proofdepth\the\secno.\the\claimno
    \or\the\secno.\the\claimno.\subclaimnumber
    \or\the\secno.\the\claimno.\subclaimnumber.\subsubclaimnumber
    \or\the\secno.\the\claimno.\subclaimnumber.\subsubclaimnumber
                                                .\subsubsubclaimnumber\fi\fi\fi}
\secno=0\claimno=0\proofdepth=0\subclaimno=0\subsubclaimno=0\subsubsubclaimno=0
\NoSectionstrue
\newbox\qedbox
\def\claimname{\ifcase\claimtype Theorem\or Lemma\or Claim\or Corollary\or
               Question\or Definition\or Remark\or Conjecture\fi}
\def\preclaimskip{\removelastskip
    \ifcase\claimtype\goodbreak\vskip 8pt plus 4pt minus 2pt
                  \or\goodbreak\vskip 6pt plus 4pt minus 1pt
                  \or\goodbreak\vskip 5pt plus 4pt minus 1pt
                  \or\goodbreak\vskip 8pt plus 4pt minus 2pt
                  \or\vskip 7pt plus 4pt minus 2pt
                  \or\vskip 7pt plus 4pt minus 2pt
                  \or\vskip 7pt plus 4pt minus 2pt
                  \or\goodbreak\vskip 8pt plus 4pt minus 2pt\fi}
\def\postclaimskip{\ifcase\claimtype         \vskip 4pt plus 2pt minus 2pt
                                          \or\vskip 3pt plus 2pt minus 2pt
                                          \or\vskip 2pt plus 2pt minus 1pt
                                          \or\vskip 4pt plus 2pt minus 2pt
                                          \or\vskip 1pt plus 2pt 
                                          \or\vskip 4pt plus 4pt 
                                          \or\vskip 3pt plus 2pt
                                          \or\vskip 4pt plus 2pt minus 2pt\fi}
\def\claimfont{\ifcase\claimtype
                  \sl\or\sl\or\sl\or\sl\or\sl\or\rm\or\rm\or\sl\fi}
\def\advancetag{\ifcase\proofdepth\advance\claimno by1
                               \or\advance\subclaimno by1
                               \or\advance\subsubclaimno by1
                               \or\advance\subsubsubclaimno by1\fi}
\def\sayclaim#1.#2 #3\par{\ifquiet\else\advancetag\fi
    \preclaimskip\setbox1=\hbox{#1}\setbox2=\hbox{#2}%
    \toks0={#1 }
    \immediate\write16{\ifdim\wd1>0pt\the\toks0
                       \else\claimname\space\fi \claimtag.}%
    \vbox{\noindent
    {\bf\ifdim\wd1=0pt \claimname\else #1\fi\ifquiet.\else\ \claimtag{\ifNoSections.\fi}\fi}%
    \enspace{\ifdim\wd2>0pt\sc #2\enspace\fi}%
    {\claimfont #3\par}}\postclaimskip\quietoff}
\def\theorem{\claimtype=0\sayclaim}
\def\lemma{\claimtype=1\sayclaim}
\def\claim{\claimtype=2\sayclaim}
\def\corollary{\claimtype=3\sayclaim}

\def\remark{\claimtype=6\sayclaim}

\def\point#1. #2\par{\item{\rm #1.}#2\par}
\def\points#1\cr\par{\medskip\vbox{\let\cr=\point\point#1\par}\par}
\def\df{\it}
\def\prooffont{}
\def\proofsize{}%\ifcase\proofdepth\or\elevenpoint\or\tenpoint\or\ninepoint\fi}
\def\proofindent{}%\ifcase\proofdepth\or\or\narrower\or\narrower\fi}
\def\proofskip{\badbreak\ifcase\claimtype    \vskip 3pt plus 2pt minus 2pt
                                          \or\vskip 2pt plus 2pt minus 2pt
                                          \or\vskip 1pt plus 2pt minus 1pt
                                          \or\vskip 3pt plus 2pt minus 2pt
                                          \or\vskip 1pt plus 2pt 
                                          \or\vskip 2pt plus 4pt 
                                          \or\vskip 1pt plus 2pt
                                          \or\vskip 3pt plus 2pt minus 2pt\fi}

\def\Goodbreak{\vskip0pt plus.5in\penalty-1000\vskip0pt plus-.5in}
\def\goodbreak{\penalty-500}
\def\badbreak{\penalty500}
\def\Badbreak{\penalty1000}
\def\proof{\message{proof}\removelastskip\Badbreak\proofskip\begingroup
  \advance\proofdepth by1
  \setbox\qedbox=\hbox{\halmos\raise2pt\hbox{\fiverm\claimname}}%
  \prooffont\proofsize\proofindent\noindent{\bf Proof: }}
\def\proofof#1:{\message{proof}\removelastskip\Badbreak\proofskip\begingroup
  \advance\proofdepth by1
  \setbox\qedbox=\hbox{\halmos\raise2pt\hbox{\fiverm#1}}%
  \prooffont\proofsize\proofindent\noindent{\bf Proof of #1: }}
\def\cite[#1]{[{\tenrm{#1}}]\message{[#1]}}
\edef\ref#1{\expandafter\global\expandafter\edef#1{\noexpand\claimtag}}
\newwrite\notes
\openout\notes=\jobname.notes
% for the following macros, see page 377 of the TeXbook
\long\def\unexpandedwrite#1#2{\def\finwrite{\write#1}%
   {\aftergroup\finwrite\aftergroup{\sanitize#2\endsanity}}}
\def\sanitize{\futurelet\next\sanswitch}
\let\stoken=\space
\def\sanswitch{\ifx\next\endsanity
  \else\ifcat\noexpand\next\stoken\aftergroup\space\let\next=\eat
   \else\ifcat\noexpand\next\bgroup\aftergroup{\let\next=\eat
    \else\ifcat\noexpand\next\egroup\aftergroup}\let\next=\eat
     \else\let\next=\copytoken\fi\fi\fi\fi \next}
\def\eat{\afterassignment\sanitize \let\next= }
\long\def\copytoken#1{\ifcat\noexpand#1\relax\aftergroup\noexpand
  \else\ifcat\noexpand#1\noexpand~\aftergroup\noexpand\fi\fi
  \aftergroup#1\sanitize}
\def\endsanity\endsanity{}

\def\note#1#2{\hbox to2in{\strut#1\quad\dotfill\quad#2}}
\def\boxit#1{\setbox4=\hbox{\kern1pt#1\kern1pt}
  \hbox{\vrule\vbox{\hrule\kern1pt\box4\kern1pt\hrule}\vrule}}
\def\halmos{\hbox{\am\char'3}} 
\def\qed#1\par{\message{.                                }\setbox1=\hbox{#1}%
  \ifdim\wd1>0pt\setbox\qedbox=\hbox{\halmos\raise2pt\hbox{\fiverm #1}}\fi
  \kern5pt\lower 2pt\hbox{\box\qedbox}\proofskip\goodbreak\endgroup}

\def\sectionsymbol{\S}
\def\k{\kappa}
\def\g{\gamma}
\def\a{\alpha}
\def\b{\beta}
\def\d{\delta}
\def\s{\sigma}
\def\t{\tau}
\def\l{\lambda}
\def\z{\zeta}
\def\I1{\mathop{\hbox{\sc i}_1}}

\def\P{{\mathchoice{\hbox{\bm P}}{\hbox{\bm P}}
         {\hbox{\tenbm P}}{\hbox{\sevenbm P}}}}
\def\Q{{\mathchoice{\hbox{\bm Q}}{\hbox{\bm Q}}
         {\hbox{\tenbm Q}}{\hbox{\sevenbm Q}}}}

\def\card#1{\left|#1\right|}

\def\dom{\mathop{\rm dom}\nolimits}

\def\unifto{\buildrel\lower 7pt\hbox{$\to$}\over\to}

\def\iso{\cong}

\def\cof{\mathop{\rm cof}\nolimits}
\def\cp{\mathop{\rm cp}\nolimits}

\def\ORD{\hbox{\sc ord}}

\def\plus{^{\scriptscriptstyle +}}

\def\in{\mathrel{\mathchoice{\raise 
1pt\hbox{$\scriptstyle\cal\char'62$}}
         {\raise 1pt\hbox{$\scriptstyle\cal\char'62$}}
         {\raise .5pt\hbox{$\scriptscriptstyle\cal\char'62$}}
         {\hbox{$\scriptscriptstyle\cal\char'62$}}}\penalty700{}}
\def\ni{\mathrel{\mathchoice{\raise 1pt\hbox{$\scriptstyle\cal\char'63$}}
                   {\raise 1pt\hbox{$\scriptstyle\cal\char'63$}}
                   {\raise .5pt\hbox{$\scriptscriptstyle\cal\char'63$}}
                   {\hbox{$\scriptscriptstyle\cal\char'63$}}}\penalty700}
\def\of{\mathrel{\mathchoice{\raise 1pt\hbox{$\scriptstyle\subseteq$}}
                   {\raise 1pt\hbox{$\scriptstyle\subseteq$}}
                   {\raise .5pt\hbox{$\scriptscriptstyle\subseteq$}}
                   {\hbox{$\scriptscriptstyle\subseteq$}}}}
\def\fo{\mathrel{\mathchoice{\raise 1pt\hbox{$\scriptstyle\supseteq$}}
                   {\raise 1pt\hbox{$\scriptstyle\supseteq$}}
                   {\raise .5pt\hbox{$\scriptscriptstyle\supseteq$}}
                   {\hbox{$\scriptscriptstyle\supseteq$}}}}
\def\notin{\mathrel{\mathchoice
  {\raise 1pt\hbox{\rlap{$\scriptstyle\;|$}$\scriptstyle\cal\char'62$}}
  {\raise 1pt\hbox{\rlap{$\scriptstyle\kern2pt 
          |$}$\scriptstyle\cal\char'62$}}
  {\raise .5pt\hbox{\rlap{$\scriptscriptstyle\, |$}$\scriptscriptstyle
      \cal\char'62$}}
  {\hbox{\rlap{$\scriptscriptstyle\, |$}$\scriptscriptstyle
     \cal\char'62$}}}%
  \penalty700}

\def\and{\mathrel{\kern1pt\&\kern1pt}}
\def\iff{\mathrel{\leftrightarrow}}

\def\Union{\bigcup}
\def\union{\cup}

\def\intersect{\cap}

\def\cross{\times}

\def\tlt{\triangleleft}
\def\nottlt{\not\!\tlt}

\def\[#1]{\left[\vphantom{\bigm|}#1\right]}
\def\<#1>{\langle\,#1\,\rangle}

\def\image{\mathbin{\hbox{\tt\char'42}}}
\def\restrict{\mathbin{\mathchoice{\hbox{\am\char'26}}{\hbox{\am\char'26}}{\hbox{\eightam\char'26}}{\hbox{\sixam\char'26}}}}
\def\force{\mathbin{\hbox{\am\char'15}}}

\def\emptyset{\mathord{\hbox{\bm\char'77}}}

\def\boolval#1{\mathopen{\lbrack\!\lbrack}\,#1\,\mathclose{\rbrack\!\rbrack}}

\def\st{\mid}
\def\seq<#1>{{\def\st{\mid\penalty650}\left<\,#1\,\right>}}

\def\set#1{\{\,#1\,\}}

\def\th{{\hbox{\fiverm th}}}

\def\lttheta{{\raise 1pt\hbox{$\scriptstyle<$}\theta}}

\def\I1{\mathop{\hbox{\sc i}_1}}
\def\ltk{{{\scriptstyle<}\k}}
\def\ltl{{{\scriptstyle<}\l}}

\def\ltd{{{\scriptstyle<}\d}}

\def\lteb{{{\scriptstyle\leq}\b}}

\def\Qdot{\dot\Q}
\def\PQ{{\P*\Qdot}}
\def\Adot{\dot A}
\def\Qdot{\dot\Q}

\def\Pforces{\force_{\P}}
\def\PQforces{\force_{\PQ}\;}

\def\Ddot{\dot D}

\def\mudot{\dot\mu}

\def\jVM{j:V\to M}

\def\Vbar{{\overline V}}

\QUIET
\def\Pkl{P_{\k}\l}
\def\Vbar{\bar V}
\def\Ddot{\dot D}
\def\Bdot{\dot B}

\def\leQ{\le_{\Q}}

\def\tddot{\dot t_\d}
\centerline{[to appear in the {\it Journal of Symbolic Logic}]}

\title Superdestructibility: A Dual to Laver's Indestructibility

\author Joel David Hamkins{\footnote\dag{\rm The first author's research
has been supported in part by the College of Staten Island and a grant from
The City University of New York PSC-CUNY Research Award Program.}}

\author Saharon Shelah{\footnote\ddag{\rm The second author's research has 
been supported by The Israel 
Science Foundation, administered by the Israel Academy of Sciences and 
Humanities. This is publication 618 in his independent numbering 
system.}}

\abstract. After small forcing, any $\ltk$-closed forcing will 
destroy the supercompactness and even the strong compactness of $\k$. 

In a delightful argument, Laver \cite[L78] proved that any supercompact 
cardinal $\k$ can be made indestructible by $\ltk$-directed closed forcing.
This indestructibility, however, is evidently not itself 
indestructible, for it is always ruined by small forcing: in \cite[H96] the 
first author recently proved that small forcing makes any 
cardinal superdestructible; that is, any further $\ltk$-closed forcing which
adds a subset to $\k$ will destroy the measurability, even the weak
compactness, of $\k$. What is more, this property holds higher up: after 
small forcing, any further $\ltk$-closed forcing which adds a subset to 
$\l$ will destroy the $\l$-supercompactness of $\k$, provided $\l$ is not 
too large (his proof needed that $\l<\aleph_{\k+\d}$, where the 
small forcing is $\ltd$-distributive). In this paper, we happily remove this 
limitation on $\l$, and show that after small forcing, the supercompactness 
of $\k$ is destroyed by any $\ltk$-closed forcing. Indeed, we will 
show that even 
the strong compactness of $\k$ is destroyed. By doing so we answer the 
questions asked at the conclusion of \cite[H96], and obtain the following
attractive complement to Laver indestructibility:

\theorem Main Theorem. After small forcing, any $\ltk$-closed forcing 
will destroy the supercompactness and even the strong compactness of $\k$.

We will provide two arguments. The first, similar to but generalizing the
Superdestruction Theorem of \cite[H96], will show that supercompactness is 
destroyed; the second,
by a different technique, will show fully that strong compactness is 
destroyed. Both arguments will rely fundamentally on the Key Lemma, below,
which was proved in \cite[H96]. Define that a set or sequence is {\df fresh} 
over $V$ when it is not in $V$ but every initial segment of it is in $V$.

\lemma Key Lemma. Assume that $\card{\P}=\b$, that $\Pforces\Qdot$ is
$\lteb$-closed, and that $\cof(\l)>\b$. Then $\P*\Qdot$ adds no fresh subsets
of $\l$, and no fresh $\l$-sequences. 

While in \cite[H96] it is proved only that no fresh sets are added, the
following simple argument shows that no fresh sequences can be added:
given a sequence in $\d^\l$, code it in the natural way with a binary
sequence of length $\d\l$, by using $\l$ many blocks of length $\d$,
each with one $1$. The binary sequence corresponds to a subset of the 
ordinal $\d\l$, which, since $\cof(\d\l)=\cof(\l)$, cannot be fresh. Thus,
the original $\l$-sequence cannot be fresh. 

Let us give now the first argument. We will use the notion of a
$\theta$-club to extend the inductive proof of the Superdestruction
Theorem \cite[H96] to all values of $\l$.

\theorem. After small forcing, any $\ltk$-closed forcing which adds a 
subset to $\l$ will destroy the $\l$-supercompactness of $\k$. 

\proof Suppose that $\card{\P}<\k$ and $\Pforces\Qdot$ is $\ltk$-closed. 
Suppose that $g*G\of\PQ$ is $V$-generic, and that $\Q=\Qdot_g$
adds a new subset $A\of\l$, with $\l$ minimal, so that $A\in V[g][G]$ but 
$A\notin V[g]$. By the closure of $\Q$, we know that $\cof(\l)\ge\k$. 
Suppose, towards a contradiction, that $\k$ is $\l$-supercompact in
$V[g][G]$. Let $\Pkl$ denote $(\Pkl)^{V[g][G]}$, which is also
$(\Pkl)^{V[g]}$. 

\lemma. Every normal fine measure on $\Pkl$ in $V[g][G]$ concentrates 
on $(\Pkl)^V$. 

\proof Let us begin with some definitions. Fix a regular cardinal $\theta$ 
such that $\card{\P}<\theta<\k$. A set $C\of\Pkl$ is 
{\df unbounded} iff for every 
$\s\in\Pkl$ there is $\t\in C$ such that $\s\of\t$. A set $D\of\Pkl$ is
{\df $\theta$-directed} iff whenever $B\of D$ and $\card{B}<\theta$ then 
there is some $\t\in D$ such that $\s\of\t$ for every $\s\in B$. 
The set $C$ is {\df $\theta$-closed} iff every $\theta$-directed 
$D\of C$ with $\card{D}<\k$ has $\union D\in C$. Finally, $C$ is a 
{\df $\theta$-club} iff $C$ is both $\theta$-closed and unbounded.

\claim. A normal fine measure on $\Pkl$ contains every $\theta$-club.

\proof Work in any model $\Vbar$. Suppose that $C$ is a $\theta$-club
in $\Pkl$ and that $\mu$ is a normal fine measure on $\Pkl$. Let 
$j:\Vbar\to M$ be the ultrapower by $\mu$. It is well known that 
$j\image\l$ is a seed for $\mu$ in the sense that 
$X\in\mu\iff j\image\l\in j(X)$ for $X\of\Pkl$. By elementarity 
$j(C)$ is a 
$\theta$-club in $M$ and $j\image C\of j(C)$. (We know $j\image C\in M$ 
because $M$ is closed under $\l^{\ltk}$ sequences in $\Vbar$.) Also, 
it is easy to check that $j\image C$ is $\theta$-directed. Thus, by 
the definition of $\theta$-club, we know $\union(j\image C)\in j(C)$. 
But $$\union(j\image C)=\Union_{\s\in C}j(\s)=
\Union_{\s\in C}(j\image\s)=j\image\l.$$
Thus, $j\image\l\in j(C)$ and so $C\in\mu$.\qed

Now let $C=(\Pkl)^V$. We will show that $C$ is a $\theta$-club in $V[g][G]$. 
First, let us show that $C$ is unbounded. If $\s\in\Pkl$ in $V[g][G]$, then
actually $\s\in V[g]$, and so $\s=\dot\s_g$ for some $\P$-name $\dot\s\in V$. 
We may assume that $\boolval{\card{\dot\s}<\check\k}=1$ and 
consequently $\s\of\set{\a\st \boolval{\a\in\dot\s}\not=0}\in C$; so 
$\s$ is covered as desired. 
To show that $C$ is $\theta$-closed, suppose in $V[g][G]$ that $D\of C$ 
has size less than $\k$ and is 
$\theta$-directed. We have to show that $\union D\in C$. It suffices
to show that $\union D\in V$ since $C=\Pkl\intersect V$. Since $\Q$ is
$\ltk$-closed, we know that $D\in V[g]$, and thus $D=\Ddot_g$ for some
name $\Ddot\in V$. In $V$ let 
$D_p=\set{\s\in C\st p\force\check\s\in\Ddot}$. It follows
that $D=\union_{p\in g}D_p$. There must be some $p\in g$ such
that $D_p$ is $\of$-cofinal in $D$; for if not, then for each $p\in g$ 
we may choose $\s_p\in D$ such that $D_p$ contains no supersets of 
$\s_p$. Since $D$ is $\theta$-directed and $\card{g}<\theta$ there is
some $\s\in D$ such that $\s_p\of\s$ for all $p\in g$. But $\s$ must
be forced into $D$ by some condition $p\in g$, so $\s\in D_p$ for 
some $p\in g$, contradicting the choice of $\s_p$. So we may fix 
some $p\in g$ such that $D_p$ is $\of$-cofinal in $D$. But in this 
case $\union D_p=\union D$ and since $D_p\in V$ we conclude $\union D\in V$. 
Thus $C$ is a $\theta$-club in $V[g][G]$, and the lemma is proved.\qed

Let us now continue with the theorem. Since $\k$ is $\l$-supercompact
in $V[g][G]$ there must be an embedding $j:V[g][G]\to M[g][j(G)]$ 
which is the ultrapower by a normal fine measure $\mu$ on $\Pkl$. 

\lemma. $P(\l)^M=P(\l)^V$. 

\proof $(\fo)$. By the previous lemma we know that $(\Pkl)^V\in\mu$ and so 
$j\image\l\in j((\Pkl)^V)=(\Pkl)^M$. Since $M$ is transitive, it follows 
that $j\image\l\in M$. And obtaining this fact was the only reason for 
proving the previous lemma. Now if $B\of\l$ and $B\in V$ then
$j(B)\in M$, and since $B$ is constructible from 
$j(B)$ and $j\image \l$ it follows that $B\in M$ as well. 

$(\of)$. Now we prove the converse. By induction we will show that 
$P(\d)^M\of V$ for all $\d\le\l$. 
Suppose that $B\of\d$ and $B\in M$ and every initial segment of $B$ 
is in $V$. By the Key Lemma it follows that $B\in V$
unless $\cof(\d)<\k$. So suppose $\cof(\d)<\k$. By the closure
of $\Q$ we know in this case that $B\in V[g]$ and so $B=\Bdot_g$ for 
some name 
$\Bdot\in V$. We may view $\Bdot$ as a function from $\d$ to the set
of antichains of $\P$. Since $\Bdot$ may be coded with a subset of
$\d$, we know $\Bdot\in M$ by the previous direction of this lemma. 
Thus, both $B$ and $\Bdot$ are in $M$ and $g$ is $M$-generic. Since 
$B=\Bdot_g$ in $M[g]$ there is in $M$ a condition $p\in g$ such that 
$p\force\Bdot=\check B$. That is,
$p$ decides every antichain of $\Bdot$ in a way that makes it 
agree with $B$. Use $p$ to decide $\Bdot$ in $V$ and conclude that 
$B\in V$. This completes the induction.\qed

Now we are nearly done. Consider again the new set $A\of\l$ such that 
$A\in V[g][G]$ but $A\notin V[g]$. Since $j$ is a $\l$-supercompact
embedding, we know $A\in M[g][j(G)]$. Since the $j(G)$ forcing is 
${<}j(\k)$-closed, we know $A\in M[g]$. Therefore $A=\Adot_g$ for 
some name $\Adot\in M$. Viewing $\Adot$ as a function from $\l$ to 
the set of antichains in $\P$, we can code $\Adot$ with a subset of
$\l$, and so by the last lemma we know $\Adot\in V$. Thus, 
$A=\Adot_g\in V[g]$, contradicting the choice of $A$.\qed

\corollary. By first adding in the usual way a generic subset to $\b$
and then to $\l$, where $\cof(\l)>\b$, one 
destroys all supercompact cardinals between $\b$ and $\l$. 

In fact, one does not even need to add them in the usual way. 
This is because the proof of the 
theorem does not really use the full $\ltk$-closure of $\Q$. Rather, 
if $\P$ has size $\b$, then we only need that $\Q$ is $\lteb$-closed
and adds no new elements of $\Pkl$. Thus, we have actually proved 
the following theorem.

\theorem. After any forcing of size $\b<\k$, any further $\lteb$-closed 
forcing which adds a subset to $\l$ but no elements to $\Pkl$ will destroy 
the $\l$-supercompactness of $\k$. 

This improvement is striking when $\b$ is small, having the 
consequence that after adding a Cohen real, any countably-closed forcing
which adds a subset to some minimal $\l$ destroys all supercompact
cardinals up to $\l$. 

Let us now give the second argument, which will improve the previous results 
with a different technique and 
establish fully that strong compactness is destroyed. 

\theorem. After small forcing, any $\ltk$-closed forcing which 
adds a $\l$-sequence will destroy the $\l$-strong compactness of $\k$.

\proof Define that a cardinal $\k$ is {\df $\l$-measurable} iff there is a 
$\k$-complete (non $\k\plus$-complete) uniform measure on $\l$. 
Necessarily $\k\leq\cof(\l)$. This notion is studied in \cite[K72]. 

\lemma. Assume that $\card{\P}<\k\leq\l$, that $\Qdot$ adds a new 
$\l$-sequence over $V^\P$, $\l$ minimal, and that $\k$ is 
$\l$-measurable in $V^{\PQ}$. Then $\PQ$ must add a fresh 
$\l$-sequence over $V$. 

\proof This lemma is the heart of the proof. 
Assume the hypotheses of the lemma. So $\PQforces\,\dot s$ is a 
$\l$-sequence of ordinals not in $V^\P$, 
and $\mudot$ is a $\k$-complete uniform measure on $\l$. Without loss of
generality, we may assume that $\Pforces\Qdot$ is a complete boolean algebra
on an ordinal. Suppose now that $g*G$ is $V$-generic for $\PQ$. Let 
$\Q=\Qdot_g$, and $s=\dot s_{g*G}$. 

In $V[g]$, let 
$T=\set{u\in\ORD^{\ltl}\st \boolval{\check u\of\dot s}^\Q\ne 0}$.
Thus, under inclusion, $T$ is a tree with $\l$ many levels, and $\Q$ adds
the $\l$-branch $s$.
For $u\in T$, let $b_u=\boolval{u\of\dot s}^\Q$.
Thus, $b_u$ is an ordinal. Let $I=\set{\<\ell(u),b_u>\st u\in T}$, where 
$\ell(u)$ denotes the length of $u$, and define 
$\<\a,b_u>\tlt\<\a',b_{u'}>$ when $\a'<\a$ and $b_u\leQ b_{u'}$. Since 
$u\supset v\iff\<\ell(u),b_u>\tlt\<\ell(v),b_v>$ it follows 
that $\<T,\supset>\iso\<I,\tlt>$, and consequently $I$ is also a 
tree, under the relation $\tlt$, with $\l$ many levels. Furthermore, 
the $\a^\th$ level of $I$ consists of pairs of the form $\<\a,\b>$. 
For $p\in\P$ let us define that
$a\tlt_p b$ when $p\force a\tlt b$. Thus, $\tlt=\union_{p\in g}\tlt_p$. 

In $V[g][G]$ let $b_\g=\langle\g,b_{s\restrict\g}\rangle$. Thus, $b_\g\in I$,
and if 
$\g<\z$ then $b_\z\tlt b_\g$ and so there is some $r\in g$ such that
$b_\z\tlt_r b_\g$. Since there are fewer than $\k$ many such $r$, 
for each $\g$ there must be an $r$ which works for $\mu$-almost 
every $\z$. But then again, since there are relatively few $r$, 
it must be that there is some $r^*\in g$ which has this 
property for $\mu$-almost every $\g$. So, fix $r^*\in g$ such that 
for $\mu$-almost every $\g$, for $\mu$-almost every $\z$, 
we have $b_\z\tlt_{r^*} b_\g$. Fix also a condition $\<p_0,q_0>\in g*G$ forcing
$r^*$ to have this property. Let 
$t=\seq<b_\g\st\g<\l\and\hbox{for $\mu$-a.e. }\z,\; b_\z\tlt_{r^*}b_\g>$. 
Thus, $t$ is a partial function from $\l$ to pairs of ordinals, and 
$\dom(t)\in\mu$. In particular, $\dom(t)$ is unbounded in $\l$. 

We will argue that $t$ is fresh over $V$. First, notice that $t\notin V[g]$ 
since in $V[g]$ knowing $t$ we could read off the branch $s$. Thus, 
$t\notin V$. 

Nevertheless, we will argue that every initial segment of $t$ is in $V$. 
Suppose $\d<\l$, and let $t_\d=t\restrict\d$. By the minimality of $\l$ 
it follows that $t_\d\in V[g]$, and so there is a $\P$-name $\tddot$
and a condition $\<p_1,q_1>\in g*G$, stronger than $\<p_0,q_0>$, 
forcing this name to work. Assume towards a 
contradiction that $t_\d\notin V$, and that this is forced by $p_1$. Then,
for each $r\in\P$ below $p_1$ we may choose $\g_r<\d$ such that $r$ does 
not decide $t(\g_r)$ (or whether $\g_r$ is in the domain of $t$).  
But, nevertheless, for each $r$ either for $\mu$-almost every $\z$, 
$b_\z\tlt_{r^*} b_{\g_r}$ or else for $\mu$-almost every $\z$, 
$b_\z\nottlt_{r^*} b_{\g_r}$ (but not both). In the first case it follows
that $t(\g_r)=b_{\g_r}$, and in the second it follows that 
$\g_r\notin\dom(t)$. Since there are relatively few $r$, by 
intersecting these sets of $\z$ we can find a single $\z$ which acts, 
with respect to the $\g_r$, exactly the way $\mu$-almost every $\z$ acts.
Fix such a $\z$. Thus, for each $r$ we have either 
$b_\z\tlt_{r^*} b_{\g_r}$, and consequently $t(\g_r)=b_{\g_r}$, or else 
$\g_r\notin\dom(t)$ (but not both). Notice 
that $\z$ and $b_\z$ are just some particular ordinals. Fix
some condition $\<p^*,q^*>$ below $\<p_1,q_1>$ forcing $\z$ and $b_\z$ to 
have the property we mention in the sentence before last. 
Now we will argue that this is a contradiction. Let $\g=\g_{p^*}$. 
There are two cases. First, it might happen that $b_\z\tlt_{r^*}\<\g,\b>$ 
for some ordinal $\b$. Such a situation can be 
observed in $V$. In this case, 
$\<p^*,q^*>$ forces $\b=b_{s\restrict\g}$ and therefore, by the assumption 
on $\z$, 
it also forces $t(\g)=\<\g,\b>$. Since $\tddot$ is a $\P$-name, it follows 
that $p^*\force\tddot(\check\g)=\<\check\g,\check\b>$, contrary to the choice 
of $\g=\g_{p^*}$.
Alternatively, in the second case, it may happen that 
$b_\z\nottlt_{r^*}\<\g,\b>$ 
for every $\b$. In this case, by the assumption on 
$\z$, it must be that $\<p^*,q^*>$ forces that $\g\notin\dom(t)$. 
Again, since $\tddot$ is a $\P$-name, it follows that 
$p^*\force\g\notin\dom(\tddot)$, contrary again to the choice of 
$\g=\g_{p^*}$. Thus, in either case we reach a contradiction, and so
we have proven that $\PQ$ must add a fresh $\l$-sequence.\qed

\lemma. If $\k\leq\cof(\l)$ and $\k$ is $\l$-strongly compact, then
$\k$ is $\l$-measurable.

\proof Let $\jVM$ be the ultrapower map witnessing that $\k$ is
$\l$-strongly compact. By our assumption on $\cof(\l)$, it follows that
$\sup j\image\l<j(\l)$. Let $\a=(\sup j\image\l)+\k$, and let $\mu$ be
the measure germinated by the seed $\a$. That is, $X\in\mu$ iff $\a\in
j(X)$. Since $\a<j(\l)$ it follows that $\mu$ is a measure on $\l$.
Since $j(\b)<\a$ for all $\b<\l$ it follows that $\mu$ is uniform.
Since $\cp(j)=\k$ it follows that $\mu$ is $\k$-complete. For $\g<\k$,
let $B_\g=\set{\b\st \g<\cof(\b)<\k}$. Since $\cof(\a)=\k$ in $M$, it
follows that $\a\in j(B_\g)$ and consequently $B_\g\in\mu$ for every
$\g<\k$.  Since $\intersect_\g B_\g=\emptyset$, it follows that $\mu$
is not $\k\plus$-complete, as desired.\qed

\remark. Ketonen \cite[K72] has proved that if $\k$ is $\l$-measurabile
for every regular $\l$ above $\k$, then $\k$ is strongly compact.  This
cannot, however, be true level-by-level, since if $\k<\l$ are both
measurable, with measures $\mu$ and $\nu$, then $\mu\cross\nu$ is a
$\k$-complete, non-$\k\plus$-complete, uniform measure on $\k\cross\l$.
Thus, in this situation, $\k$ will be $\l$-measurable, even when it may
not be even $\k\plus$-strongly compact.  But the previous lemma
establishes that the direction we need does indeed hold
level-by-level.

Let us now finish the proof of the theorem. Suppose that $V[g][G]$ is a
forcing extension by $\PQ$, where $\card{\P}<\k$ and $\Q$ is
$\ltk$-closed.  Let $\l$ be least such that $\Q$ adds a new
$\l$-sequence not in $V[g]$. Necessarily, $\k\leq\l$ and $\l$ is
regular. By the Key Lemma $V[g][G]$ has no $\l$-sequences which are
fresh over $V$. Thus, by the first lemma $\k$ is not $\l$-measurable in
$V[g][G]$. Therefore, by the second lemma, $\k$ is not $\l$-strongly
compact in $V[g][G]$.\qed

So the proof actually establishes that after small forcing of size
$\b<\k$, any $\lteb$-closed forcing which adds a new $\l$-sequence for
some minimal $\l$, with $\l\ge\k$, will destroy the $\l$-measurability
of $\k$. This subtlety about adding a $\l$-sequence as opposed to a
{\it subset} of $\l$ has the following intriguing consequence, which is
connected with the possibilities of changing the cofinalities of very
large cardinals.

\corollary. Suppose that $\k$ is $\l$-measurable. Then after forcing
with $\P$ of size $\b<\k$, any $\lteb$-closed $\Q$ which adds a
$\l$-sequence, but no shorter sequences, must necessarily add subsets to $\l$.

\proof Such forcing will destroy the $\l$-measurability of $\k$. Hence,
it must add subsets to $\l$.\qed

\section Bibliography

\tenpoint
\nopagenumbers
\parindent=0pt
\newbox\Article
\newbox\Journal
\newbox\Author
\newbox\Vol
\newbox\No
\newbox\Year
\newbox\Page
\newbox\Book
\newbox\Publisher
\newbox\Pubaddr
\newbox\Key
\newbox\Editor
\newbox\Comment
\newbox\Note
\def\entry#1#2\par{\item{#1\quad}\hskip-1.1em#2\par}
\def\article#1{\setbox\Article=\hbox{\sl #1, }}
\def\journal#1{\setbox\Journal=\hbox{\rm #1 }}
\def\author#1{\setbox\Author=\hbox{\sc #1, }}
\def\vol#1{\setbox\Vol=\hbox{\bf #1 }}
\def\no#1{\setbox\No=\hbox{no. #1 }}
\def\year#1{\setbox\Year=\hbox{\rm({\oldstyle #1}) }}
\def\page#1{\setbox\Page=\hbox{\rm p. #1 }}
\def\book#1{\setbox\Book=\hbox{\it #1, }}
\def\publisher#1{\setbox\Publisher=\hbox{\rm #1, }}
\def\pubaddr#1{\setbox\Pubaddr=\hbox{\rm #1, }}
\def\key#1{\setbox\Key=\hbox{#1}}
\def\editor#1{\setbox\Editor=\hbox{\rm(#1, Ed.) }}
\def\comment#1{\setbox\Comment=\hbox{\rm #1}}
\def\note#1{\setbox\Note=\hbox{\rm #1 }}
\def\ref#1\par{\smallskip{#1
  \entry{\ifhbox\Key\unhbox\Key\else[\ ]\fi}%
  \unhbox\Author\unhbox\Note
  \ifhbox\Book \unhbox\Book\unhbox\Publisher\unhbox\Pubaddr
               \unhbox\Editor\unhbox\Page\unhbox\Year\unhbox\Comment
  \else \unhbox\Article\unhbox\Journal\unhbox\Vol\unhbox\No\unhbox\Editor
        \unhbox\Page\unhbox\Year\unhbox\Comment\fi\par}}

\ref
\author{Joel David Hamkins}
\article{Small Forcing Makes Any Cardinal Superdestructible}
\journal{Journal of Symbolic Logic}
\comment{(in press)}
\key{[H96]}

\ref
\author{Jussi Ketonen}
\article{Strong Compactness and Other Cardinal Sins}
\journal{Annals of Mathematical Logic}
\vol{5}
\year{1972}
\page{47-76}
\key{[K72]}

\ref
\author{Richard Laver}
\article{Making the Supercompactness of $\kappa$ Indestructible Under 
 $\kappa$-Directed Closed Forcing}
\journal{Isreal Journal Math}
\vol{29}
\year{1978}
\page{385-388}
\key{[L78]}

\bye